\documentclass{gtart_h}

\def\ifplaintex{\expandafter\ifx\csname documentclass\endcsname\relax}

\ifplaintex 
\else
\fi

\def\gt{{\mathsurround=0pt\it $\cal G\mskip-2mu$eometry \&\ 
$\cal T\!\!$opology}}        

\def\gtm{{\mathsurround=0pt\it $\cal G\mskip-2mu$eometry \&\ 
$\cal T\!\!$opology $\cal M\mskip-1mu$onographs}}    

\def\gtp{{\mathsurround=0pt\it $\cal G\mskip-2mu$eometry \&\ 
$\cal T\!\!$opology $\cal P\!$ublications}}  

\def\recd{{\small Received:\qua\receiveddate\ifx\reviseddate\relax
\else\qquad Revised:\qua\reviseddate\fi\par}} 

\def\volumenumber#1{\def\thevolumenumber{#1}}
\def\volumeyear#1{\def\thevolumeyear{#1}}
\def\volumename#1{\def\thevolumename{#1}}
\def\papernumber#1{\def\thepapernumber{#1}}
\def\pagenumbers#1#2{\def\startpage{#1}\def\finishpage{#2}}
\def\published#1{\def\publishdate{#1}}
\def\received#1{\def\receiveddate{#1}}
\def\revised#1{\def\reviseddate{#1}}
\def\accepted#1{\def\accepteddate{#1}}

\long\def\asciiabstract#1{\long\def\theasciiabstract{#1}}

\let\\\par
\let\thevolumenumber\relax\let\thepapernumber\relax
\let\thevolumeyear\relax\let\startpage\relax
\let\finishpage\relax\let\publishdate\relax\let\receiveddate\relax
\let\reviseddate\relax\let\accepteddate\relax\let\theasciititle\relax
\let\theasciiauthors\relax
\let\theasciiabstract\relax

\let\theerratum\relax\let\theasciiemail\relax
\let\theshortauthors\relax\let\theshorttitle\relax

\def\startpage{1}\def\finishpage{15}\def\thepapernumber{77}

\volumenumber{2}
\volumename{Proceedings of the Kirbyfest}
\volumeyear{1999}

\long\def\maketitlep{   

\count0=\startpage

\gtm\nl        
{\small Volume \thevolumenumber: \thevolumename\nl 
\ifx\theerratum\relax\else Erratum \erratumnumber\nl\fi
Pages \startpage--\finishpage\nl}

\vglue 0.1truein   

{\parskip=0pt\leftskip 0pt plus 1fil\def\\{\par\smallskip}{\ifplaintex\large
\else\Large\fi\bf\thetitle}\par\medskip}   
\vglue 0.05truein 

{\parskip=0pt\leftskip 0pt plus 1fil\def\\{\par}{\sc\theauthors}
\par\medskip}%
 
\vglue 0.03truein 

{\small\leftskip 25pt\rightskip 25pt{\bf Abstract}\stdspace\theabstract

{\bf AMS Classification}\stdspace\theprimaryclass
\ifx\thesecondaryclass\relax\else; \thesecondaryclass\fi\par
{\bf Keywords}\stdspace \thekeywords\par}\vglue 7pt

}   

\font\phead=cmsl9 scaled 950
\font=cmsl9 scaled 1050
\font\pnum=cmbx10 scaled 913
\font\lnum=cmbx10 
\font\pfoot=cmsl9 scaled 950
\font=cmsl9 scaled 1050
\ifplaintex
\headline{\vbox to 0pt{\vskip -4.5mm\line{\small\phead\ifnum
\count0=\startpage ISSN 1464-8997 (on line)
1464-8989 (printed) \hfill {\pnum\folio}\else\ifodd\count0\def\\{ }%
\ifx\theshorttitle\relax\thetitle\else\theshorttitle\fi\hfill{\pnum\folio}
\else\def\\{ and }{\pnum\folio}\hfill\ifx\theshortauthors\relax\theauthors
\else\theshortauthors\fi\fi\fi}\vss}}
\footline{\vbox to 0pt{\vglue 0mm\line{\small\pfoot\ifnum\count0=\startpage
Published \publishdate:\qua\copyright\ \gtp\hfill\else
\gtm, Volume \thevolumenumber\ (\thevolumeyear)\hfill\fi}\vss
}}
\else
\makeatletter
\def\@oddhead{{\small\lhead\ifnum\count0=\startpage ISSN 1464-8997 (on line)
1464-8989 (printed) \hfill {\lnum\number\count0}\else\ifodd\count0
\def\\{ }\ifx\theshorttitle\relax \thetitle \else\theshorttitle\fi\hfill
{\lnum\number\count0}\else\def\\{ and }{\lnum\number\count0}
\hfill\ifx\theshortauthors\relax 
\theauthors\else\theshortauthors\fi\fi\fi}}\def\@evenhead{@oddhead}
\def\@oddfoot{\small\lfoot\ifnum\count0=\startpage Published \publishdate:\qua\copyright\ \gtp\hfill\else
\gtm, Volume \thevolumenumber\ (\thevolumeyear)\hfill\fi}
\def\@evenfoot{@oddfoot}
\makeatother
\fi

\let\maketitlepage\maketitlep
\let\makeshorttitle\maketitlepage
\let\maketitle\maketitlepage

\newwrite\gtoutfile
\long\gdef\makeheadfile{  
{\def\\{, }\def\s{ }
\immediate\openout\gtoutfile head.xxx
\immediate\write\gtoutfile{Proxy-for: \ifx\theasciiauthors\relax
\theauthors\else\theasciiauthors\fi\s<\ifx\theasciiemail\relax\theemail\else\theasciiemail\fi>}
\immediate\write\gtoutfile{\noexpand\\}
\immediate\write\gtoutfile{Authors: \ifx\theasciiauthors\relax
\theauthors\else\theasciiauthors\fi}
{\def\\{ }\immediate\write\gtoutfile{Title: \ifx\theasciititle\relax
\thetitle\else\theasciititle\fi}}
\immediate\write\gtoutfile{Subj-class: GT or SG, GR etc}
\immediate\write\gtoutfile{MSC-class: \theprimaryclass\ifx\thesecondaryclass\relax\else, \thesecondaryclass\fi}
\immediate\write\gtoutfile{Journal-ref: Geom. Topol. Monogr. \thevolumenumber\s
(\thevolumeyear) \startpage-\finishpage}
\immediate\write\gtoutfile{Comments: Published by Geometry and Topology Monographs at}
\immediate\write\gtoutfile{\s\s\s  http://www.maths.warwick.ac.uk/gt/GTMon\thevolumenumber/paper\thepapernumber.abs.html}
\immediate\write\gtoutfile{\noexpand\\}
\immediate\write\gtoutfile{}
\ifx\theasciiabstract\relax
\immediate\write\gtoutfile{\theabstract}\else
\immediate\write\gtoutfile{\theasciiabstract}\fi
\immediate\write\gtoutfile{}
\immediate\write\gtoutfile{\noexpand\\}
\immediate\write\gtoutfile{}
\immediate\closeout\gtoutfile}}  

\def\maketitlepage{\maketitlep\makeheadfile}
\let\makeshorttitle\maketitlepage
\let\maketitle\maketitlepage

\volumenumber{7}
\volumename{Proceedings of the Casson Fest}
\volumeyear{2004}
\papernumber{11}
\pagenumbers{267}{290}
\received{4 June 2004}
\revised{2 August 2004}
\accepted{20 July 2004}
\published{20 September 2004}

\usepackage{amsmath,amsfonts,amssymb}

\newtheorem{prop}{Proposition}[section]
\newtheorem{thm}[prop]{Theorem}

\newtheorem{lem}[prop]{Lemma}

\newtheorem{adden}[prop]{Addendum}
\theoremstyle{definition}
\newtheorem{de}[prop]{Definition}

\newtheorem{prince}[prop]{Principle}
\theoremstyle{remark}
\gdef\mbox#1{\leavevmode\hbox{#1}}
\def\without{\kern.4pt\mbox{%
   \vrule width 1ex height.08ex depth0ex\kern-0.08ex
      \vrule width0.08ex height1.5ex depth0ex}\kern.4pt}
      \def\contract{{\ \without\ }}
\def\D{{\mathcal D}}
\def\J{{\mathcal J}}
\def\cl{\text{\rm cl}}
\def\GL{\text{\rm GL}}
\def\dR{\text{\rm dR}}
\def\rU{\text{\rm U}}

\def\ep{\varepsilon}
\def\cee{{\mathbb C}}
\def\complex{{\mathbb C}}

\def\CP{{\mathbb C\mathbb P}}
\def\real{{\mathbb R}}
\def\zed{{\mathbb Z}}
\def\Aut{\operatorname{Aut}}
\def\codim{\operatorname{codim}}
\def\id{\operatorname{id}}
\def\Im{\operatorname{Im}}
\def\Int{\operatorname{int}}
\def\rel{\operatorname{rel}}
\def\std{\operatorname{std}}
\def\supp{\operatorname{supp}}
\def\ostd{{\omega_{\std}}}

\begin{document}

\title{Symplectic structures from Lefschetz pencils\\in high dimensions}
\author{Robert E Gompf}
\address{Department of Mathematics, The University of Texas at Austin\\1
University Station C1200, Austin, TX 78712--0257, USA}
\email{gompf@math.utexas.edu}
\keywords{Linear system, vanishing cycle, monodromy}
\subjclass{57R17}

\begin{abstract}
A symplectic structure is canonically constructed on any manifold endowed 
with a topological  linear $k$--system whose fibers carry suitable 
symplectic data. 
As a consequence, the classification theory for Lefschetz pencils in the 
context of symplectic topology is analogous to the corresponding theory 
arising in differential topology. 
\end{abstract}
\asciiabstract{%
A symplectic structure is canonically constructed on any manifold
endowed with a topological linear k-system whose fibers carry
suitable symplectic data.  As a consequence, the classification theory
for Lefschetz pencils in the context of symplectic topology is
analogous to the corresponding theory arising in differential
topology.}

\maketitle

\section{Introduction}

There is a classical dichotomy between flexible, topological objects such as 
smooth manifolds, and rigid, geometric objects such as complex 
algebraic varieties. 
Symplectic manifolds lie somewhere between these two extremes, raising  the 
question of whether they should be considered as fundamentally topological 
or geometric. 
One approach to this question can be traced back to 
Lefschetz, who attempted  to bridge the gap between topology and algebraic 
geometry by introducing topological (fibrationlike) structures now called 
{\em Lefschetz pencils} on any algebraic variety. 
These structures and more general {\em linear systems} 
can also be defined in the setting of differential topology, 
where they can be found on many manifolds that do not admit algebraic 
structures, and provide deep information about the topology of the underlying
manifolds. 
It is now becoming apparent that the appropriate context for studying linear 
systems is not algebraic geometry, but a larger context that includes all 
symplectic manifolds. 
Every closed symplectic manifold (up to deformation) admits linear 1--systems
(Lefschetz pencils)
\cite{D} and 2--systems \cite{A}, and it seems reasonable to expect 
linear $k$--systems for all $k$.
Conversely, linear $(n-1)$--systems on smooth $2n$--manifolds determine
symplectic structures \cite{G1}. 
In this paper, we show that for any $k$, a linear $k$--system, endowed with 
suitable symplectic data on the fibers, determines  a symplectic structure 
on the underlying manifold (Theorem~\ref{thm:taming}). 
We then apply this to the study of  Lefschetz pencils, 
to provide a framework in which symplectic 
structures appear much more topological than algebrogeometric. 
While Lefschetz pencils in the algebrogeometric world carry delicate 
algebraic structure, topological Lefschetz pencils have a classification 
theory expressed entirely in terms of embedded spheres and a diffeomorphism 
group of the fiber. 
The main conclusion of this article 
(Theorem~\ref{thm:omega}) is that symplectic Lefschetz pencils 
have an analogous classification theory in terms of Lagrangian spheres and 
a symplectomorphism group of the fiber. 
That is, the subtleties of symplectic geometry do not interfere with a 
topological approach to classification. 

To construct a prototypical linear $k$--system on a smooth algebraic 
variety $X\subset\CP^N$ of complex dimension~$n$, simply choose a linear
subspace $A\subset \CP^N$ of codimension $k+1$, with $A$ transverse to $X$. 
The {\em base locus\/} $B= X\cap A$ is a complex submanifold of $X$ with 
codimension $k+1$. 
The subspace $A\subset \CP^N$ lifts to a codimension--$(k+1)$ linear 
subspace $\widetilde{A}\subset \cee^{N+1}$, and projection to
$\cee^{N+1}/\widetilde{A} \cong \cee^{k+1}$ descends to a holomorphic
map $\CP^N -A\to\CP^k$ whose restriction will be denoted $f\co X-B\to \CP^k$.
The {\em fibers\/} $F_y = f^{-1} (y) \cup B$ of this linear $k$--system 
are the intersections of $X$ with the codimension--$k$ linear subspaces of 
$\CP^N$ containing $A$. 
The transversality hypothesis guarantees that $B\subset X$ has a tubular 
neighborhood $V$ with a complex vector bundle structure $\pi\co V\to B$ such 
that $f$ restricts to projectivization $\cee^{k+1} - \{0\} \to\CP^k$ 
(up to action by $\GL(k+1,\cee)$) on each fiber.

To generalize this structure to a smooth $2n$--manifold $X$, we first need to 
relax the holomorphicity conditions. 
Recall that an {\em almost-complex structure} $J\co TX\to TX$ on $X$ is a 
complex vector bundle structure on the tangent bundle (with each 
$J_x\co T_x X\to T_x X$ representing multiplication by $i$). 
This is much weaker than a holomorphic structure on $X$. 
For our purposes, it is sufficient to assume $J$ is continuous (rather 
than smooth). 
We impose such a structure on $X$, but rather than requiring $f\co
X-B\to \CP^k$ to be $J$--holomorphic (complex linear on each tangent
space), it suffices to impose a weaker condition. 
Let $\ostd$ denote the standard (K\"ahler) symplectic structure on $\CP^k$, 
normalized so that $\int_{\CP^1} \ostd =1$. 
(Recall that a symplectic structure is a closed 2--form that is 
nondegenerate as a bilinear form on each tangent space.) 
We require $J$ on $X$ to be $(\ostd,f)$--tame in the following sense:

\begin{de}\label{de:1.1}
\cite{G1}\qua 
For a $C^1$ map $f\co X\to Y$ and a 2--form $\omega$ on $Y$, an almost-complex 
structure $J$ on $X$ is {\em $(\omega,f)$--tame} if $f^* \omega (v,Jv)>0$ 
for all $v\in TX - \ker df$.
\end{de}

\noindent
In the special case $f=\id_X$, this reduces to the standard notion of $J$ 
being $\omega$--tame. 
In that case, 
imposing the additional condition that $\omega (Jv,Jw)= \omega(v,w)$ 
for all $x\in X$ and $v,w\in T_xX$ gives the notion of $\omega$--compatibility. 
For example, the standard complex structure on $\CP^k$ is 
$\ostd$--compatible, so the standard complex structure on our algebraic 
prototype $X\subset \CP^N$ is $(\ostd,f)$--tame for 
$f\co X-B\to \CP^k$ as above. 
For $f= \id_X$, 
the $\omega$--tameness condition (unlike $\omega$--compatibility) is
open, ie
preserved under small perturbations of $\omega$ and $J$, and a closed 
$\omega$ taming some $J$ is automatically symplectic (since it is 
nondegenerate: every nonzero $v\in TX$ pairs nontrivially with something, 
namely $Jv$). 
Such pairs $\omega$  and $J$ determine the same orientation on $X$.
In general, 
the $(\omega,f)$--tameness condition is preserved under taking convex 
combinations of forms $\omega$ (for fixed $J,f$). 
If $J$ is $(\omega,f)$--tame, then each $\ker df_x \subset T_xX$ is a 
$J$--complex subspace (characterized as those $v\in T_xX$ for which 
$f^*\omega (v,Jv)=0$), so away from critical points each $f^{-1}(y)$ is 
a $J$--complex submanifold of $X$.

We can now define linear systems on smooth manifolds: 

\begin{de}\label{def:1.2}
For $k\ge1$, a {\em linear $k$--system} $(f,J)$ on a smooth, closed
$2n$--manifold $X$ is a closed, codimension--$2(k+1)$ submanifold
$B\subset X$, a smooth $f\co X-B\to \CP^k$, and a continuous
almost-complex structure $J$ on $X$ with $J|X-B$ $(\omega_{\std},f)$--tame,
such that $B$ admits a neighborhood $V$ with a (smooth, correctly oriented) 
complex vector bundle structure $\pi\co V\to B$ for which $f$ is 
projectivization on each fiber.
\end{de}

For each $y\in \CP^k$, the fiber $F_y = f^{-1}(y)\cup B$ is a
closed subset of $X$ whose intersection with $V$ is a smooth,
codimension--$2k$ submanifold.
$F_y$ is a $J$--holomorphic submanifold away from the critical points of $f$,
since $J$ is $(\omega_{\std},f)$--tame on $X-B$ and continuous at $B$.
The complex orientation of $F_y$ agrees with the preimage orientation 
induced from the complex orientations of $X$ and $\CP^k$.
The base locus $B= F_y \cap F_{y'}$ ($y'\ne y\in \CP^k$)
is $J$--holomorphic.
The complex orientation of $B$, which in the transverse case $k=1$ is also 
the intersection orientation of $F_y\cap F'_y$, 
determines the ``correct'' orientation for the
fibers of $\pi$.
Later (Lemma~\ref{nu})
we will verify that the complex bundle structure on $V$ can be assumed 
to come from $J$ on $TX|B$ by the Tubular Neighborhood Theorem.

Our first goal is to construct symplectic structures using linear $k$--systems.
This was already achieved in \cite{G1} for {\em hyperpencils}, which are 
linear $(n-1)$--systems endowed with some additional structure taken from 
the algebraic prototype. 
It was shown that every hyperpencil determines a unique symplectic form up 
to isotopy. 
(Symplectic forms $\omega_0$ and $\omega_1$ on $X$ are {\em isotopic} if 
there is a diffeomorphism $\psi\co X\to X$ isotopic to $\id_X$ with 
$\psi^* \omega_0 = \omega_1$.)
The proof crucially used the fact that fibers of linear $(n-1)$--systems are 
oriented surfaces (away from the critical points) --- note that by Moser's 
Theorem \cite{M} every closed,  connected, oriented surface admits a unique 
symplectic form (ie area form) up to isotopy and scale. 
For $k< n-1$, the fibers will have higher dimension, so symplectic forms 
on them need neither exist nor be unique, and we must hypothesize existence 
and some compatibility of symplectic structures on the fibers. 
Similarly, almost-complex structures exist essentially uniquely on 
oriented surfaces, so the required almost-complex structure on a hyperpencil 
can be essentially uniquely constructed, given only a local existence 
hypothesis at the critical points. 
For higher dimensional fibers, there seems to be no analogous procedure, 
requiring us to include a global almost-complex structure in the 
defining data of a linear $k$--system. 
(Consider the projection $S^2\times S^4 \to S^2$ which can be made holomorphic
locally, but whose fibers admit no almost-complex structure.) 
The main result for constructing symplectic forms on linear $k$--systems 
is Theorem~\ref{thm:taming}. 
The statement is rather technical, but can be informally summed up as follows:

\begin{prince}\label{A}
For a linear $k$--system $(f,J)$ on $X$, suppose that the fibers admit 
$J$--taming 
symplectic structures (suitably interpreted at the critical points),
and that these can be chosen to fit together suitably 
along $B$ and in cohomology. 
Then $(f,J)$ determines an isotopy class of symplectic forms on $X$.
\end{prince}

The isotopy class of forms can be explicitly characterized 
(Addenda~\ref{unique} and \ref{Jfixed}). 

Our main application concerns Lefschetz pencils on smooth manifolds. 
These are structures obtained by generalizing the generic algebraic 
prototype of linear 1--systems. 

\begin{de}\label{de:pencil}
A {\em Lefschetz pencil} on a smooth, closed, oriented $2n$--manifold $X$
is a closed, codimension--4 submanifold $B\subset X$ and a smooth 
$f\co X-B\to \CP^1$ such that
\begin{itemize}
\item[(1)] $B$ admits a neighborhood $V$ with a 
(smooth, correctly oriented) complex
vector bundle structure $\pi\co V\to B$ for which $f$ is projectivization on
each fiber,
\item[(2)] for each critical point $x$ of $f$, there are
orientation-preserving coordinate charts about $x$ and $f(x)$ 
(into $\cee^n$ and $\cee$, respectively) in which
$f$ is given by $f(z_1,\ldots,z_n) = \sum_{i=1}^n z_i^2$, and
\item[(3)] $f$ is 1--1 on the critical set $K\subset X$.
\end{itemize}
\end{de}

Condition (2) implies $K$ is finite, so (3) can always be achieved by
a perturbation of $f$.
A Lefschetz pencil, together with an $(\ostd,f)$--tame almost-complex
structure $J$, is a linear 1--system (although the latter can have more
complicated critical points).
Such a Lefschetz 1--system can be constructed as before 
on any smooth algebraic variety
by using a suitably generic linear subspace
$A\cong \CP^{N-2} \subset \CP^N$.
On the other hand, projection $S^2 \times S^4\to S^2 = \CP^1$ 
gives a (trivial) Lefschetz pencil admitting no such $J$.

The topology of Lefschetz pencils is understood at the most basic level, 
eg \cite{L} or (in dimension~4) \cite{GS}. 
We first consider the case with $B=\emptyset$, or {\em Lefschetz fibrations}
$f\co X^{2n}\to S^2$. 
Choose a collection $A= \bigcup  A_j \subset S^2$ of embedded arcs with 
disjoint interiors, connecting the critical values to a fixed regular value 
$y_0 \in S^2$. 
Over a sufficiently small disk $D\subset S^2$ containing $y_0$, we see 
the trivial bundle $D\times F_{y_0}\to D$. 
Expanding $D$ to include an arc $A_j$ adds an $n$--handle along an 
$(n-1)$--sphere lying in a fiber. 
Thus, if we expand $D$ to include $A$, the result is specified by a 
cyclically ordered collection of {\em vanishing cycles}, ie embeddings 
$S^{n-1}\to F_{y_0}$ with suitable normal data.
However, this ordered collection depends on our choice of $A$. 
Any change in $A$ can be realized by a sequence  of {\em Hurwitz moves}, 
moving some arc $A_j$ past its neighbor $A_{j\pm1}$. 
The effect of a Hurwitz move on the ordered collection of vanishing cycles 
can be easily described using the monodromy of the fibration around 
$A_{j\pm1}$, which is an explicitly understood element of $\pi_0$ of the
diffeomorphism group $\D$ of $F_{y_0}$. 
(See Section~3.) 
Over the remaining disk $S^2 - \text{int }D$, we again have a trivial bundle, 
so the product of the monodromies of the vanishing cycles must be trivial, 
and then the Lefschetz fibrations extending fixed data over $D$ are 
classified by $\pi_1(\D)$. 
The correspondence with $\pi_1(\D)$ is determined by fixing an arc from 
$y_0$ to $\partial D$ (avoiding $A$) and a trivialization of $f$ over 
$\partial D$. 
Hurwitz moves involving the new arc will induce additional equivalences. 
For the case $B\ne\emptyset$, we blow up $B$ to obtain a Lefschetz fibration, 
then apply the previous analysis.
However, extra care is required to preserve the blown up base locus and its 
normal bundle. 
We must take $\D$ to be the group  of diffeomorphisms of $F$ fixing $B$ 
and its normal bundle, and the product of monodromies will now be a 
nontrivial normal twist $\delta$ around $B$. 
We state the result carefully as Proposition~\ref{prop:fixedA}. 
For now, we sum up the discussion as follows:

\begin{prince}\label{B}
To classify Lefschetz pencils with a fixed fiber and base locus, first 
classify, up to Hurwitz moves, cyclically ordered collections of vanishing 
cycles for which the product of monodromies is $\delta \in\pi_0(\D)$. 
For any fixed choice of arcs and vanishing cycles, 
the resulting Lefschetz pencils are classified by $\pi_1(\D)$. 
The final classification results from modding out the effects of Hurwitz moves 
on the last fiber. 
(One may also choose to mod out by self-diffeomorphisms of the fiber 
$(F_{y_0},B)$.)
\end{prince}

Of course, this is an extremely difficult problem in general, but at least 
we know where to start. 

If $X$ is given a symplectic structure $\omega$ that is symplectic on 
the fibers, then the above description can be refined. 
The vanishing cycles will be {\em Lagrangian} spheres (ie $\omega$ 
restricts to $0$ on them), and the monodromies will be symplectomorphisms
(diffeomorphisms preserving $\omega$) \cite{Ar,S1,S2}. 
The discussion of arcs and Hurwitz moves proceeds as before, where $\D$ 
is replaced by a suitable group $\D_{\omega_F}$ 
of symplectomorphisms of the fiber. 
However, symplectic forms are {\em a~priori} global analytic objects
(satisfying the partial differential equation $d\omega=0$), so for
symplectic forms on $X$
compatible with a given Lefschetz pencil, one might expect both the 
existence and uniqueness questions to involve delicate analytic invariants. 
Our main result (Theorem~\ref{thm:omega}) is that no such difficulties arise, 
provided that we choose our definitions with suitable care, for example 
requiring $[\omega]\in H_{\dR}^2 (X)$ to be Poincar\'e dual to the fibers 
(as is the case for Donaldson's pencils \cite{D}). 
We obtain: 

\begin{prince}\label{C}
The classification of (suitably defined) symplectic Lefschetz pencils is 
purely topological, ie analogous to that of Principle~\ref{B}. 
More precisely, for a suitable symplectic manifold pair $(F,B)$, let 
$i_*$ denote the $\pi_1$--homomorphism induced by inclusion 
$\D_{\omega_F}\subset \D$. 
Then for fixed (suitably symplectic) data over $D$ as preceding 
Principle~\ref{B}, a given Lefschetz pencil admits a suitably compatible 
symplectic structure if and only if it is classified by an element of 
$\Im i_*$. 
Then such structures are classified up to suitable isotopy 
by $\pi_2(\D/\D_{\omega_F})$, and by  $\ker i_*$  if symplectomorphisms 
preserving $f$ and fixing $f^{-1}(D)$ are also allowed. 
\end{prince}

This is the same sort of topological classification one obtains for 
extending bundle structures over a 2--cell: 
Given groups $H\subset G$, a space $Y\cup 2$--cell, and a fixed $H$--bundle 
over $Y$ (on which we do not allow automorphisms), $G$--bundle and 
$H$--bundle extensions (if they exist) are classified by 
$\pi_1(G)\cong \pi_2 (BG)$ and $\pi_1 (H)\cong \pi_2 (BH)$, respectively. 
Inclusion $i\co H\to G$  induces an exact sequence 
$$\pi_2 (G/H) \xrightarrow{\ \partial_*\ } \pi_1(H)  \xrightarrow{\ i_*\ } 
\pi_1 (G)$$
with $i_*$ corresponding to the forgetful map from $H$--structures to 
$G$--structures. 
Thus $\Im i_*$ classifies $G$--extensions admitting $H$--reductions, and 
$\ker i_* = \Im \partial_*$ classifies $H$--reductions of a fixed 
$G$--extension as abstract $H$--extensions. 
However, different $H$--reductions can be abstractly $H$--isomorphic via 
a $G$--bundle automorphism supported over the 2--cell, and if we disallow 
such equivalences, $H$--reductions of a fixed $G$--extension are classified 
by $\pi_2 (G/H)$.

\section{Linear systems}

In this section, we show how to construct symplectic structures from linear 
systems with suitable symplectic data along the fibers 
(Principle~\ref{A}). 
Our construction is modeled on the corresponding method for hyperpencils 
\cite[Theorem 2.11]{G1}, but is complicated by the fact that the base locus 
need no longer be $0$--dimensional. 
We must first gain more control of the normal data along $B$. 
Given a linear $k$--system $(f,J)$ on $X$ as in Definition~\ref{def:1.2}, 
let $\nu \to B$ be any $J$--complex subbundle of $TX|B$ complementary to $TB$. 
(This exists since $B$ is a $J$--holomorphic submanifold of $X$.) 
Then the bundle structure $\pi\co V\to B$ guaranteed on a neighborhood of $B$ 
(by Definition~\ref{def:1.2}) can be arranged 
(after precomposing $\pi$ with an isotopy 
preserving $f$) to have its fibers tangent to $\nu$ along $B$.

\begin{lem}\label{nu}
For $\nu$ and $\pi$ as above, the complex bundle structure on $\pi$ 
(given by Definition~\ref{def:1.2}) restricts to $J$ on $\nu$. 
\end{lem}

\begin{proof} 
Near $B$, extend $TB$ to a $J$--complex subbundle $H$ of $TX$ complementary 
to the fibers of $\pi$ and tangent to the fibers $F_y$ of $f$. 
Then $J$ induces a complex structure near $B$ on $TX/H$. 
The latter bundle is canonically $\real$--isomorphic to the bundle of tangent 
spaces to the fibers of $\pi$; let $J'$ denote the resulting almost-complex 
structure on the fibers of $\pi$. 
Clearly, $J'=J$ on $\nu$, so it suffices to show that $J'$ agrees with 
the complex structure of $\pi$ on $\nu$. 
This follows immediately from \cite[Lemma 4.4(b)]{G1}, which is restated below.
(Note that for $x\notin B$, $H_x$ lies in $\ker df_x$, so $J'$ is 
$(\ostd,f)$--tame at $x$ since $J$ is.)
\end{proof}

\begin{lem}\label{lem:lines}
\cite{G1}\qua
If $f\co\cee^n-\{0\}\to \CP^{n-1}$ denotes
projectivization, $n\ge 2$, and $J$ is a continuous (positively oriented) 
almost-complex 
structure on a neighborhood $W$ of $0$ in $\cee^n$,  with $J|W-\{0\}$
$(\omega_{\std},f)$--tame, then $J|T_0\cee^n$ is the standard complex
structure.
\end{lem}

The main idea of the proof is that $J|T_0\cee^n$ has the same complex
lines as the standard structure (since the complex lines of $\cee^n$ are
$J$--complex by $(\ostd,f)$--tameness), and a linear complex structure
is determined by its complex lines.

We can now state the main theorem of this section. 
By Lemma~\ref{nu}, the canonical identification of the vector bundle 
$\pi\co V\to B$ with the normal bundle $\nu$ of its $0$--section is 
a $J$--complex isomorphism. 
This complex bundle is projectively trivialized by $f$ 
(in Definition~\ref{def:1.2}), so we can 
reduce the structure group of $\nu$ to $\rU(1)$  (acting diagonally
on $\cee^{k+1}$) by choosing a suitable Hermitian structure on $\nu$.
This Hermitian structure is canonically determined up to a positive 
scalar function. 
Let $h$ denote the hyperplane class in $H_{\dR}^2 (\CP^k)$ dual 
to $[\CP^{k-1}]$, and let 
$c_f\in H_{\dR}^2 (X)$ correspond to 
$f^* h\in H_{\dR}^2 (X-B)$ under the obvious isomorphism. 
(Recall $\codim B\ge 4$.) 

\begin{thm}\label{thm:taming}
Let $(f,J)$ be a linear $k$--system on $X$.
Choose a $J$--complex subbundle $\nu\subset TX|B$ complementary to $TB$, 
and a Hermitian form on $\nu$ as above. 
Suppose there is a symplectic form $\omega_B$ on $B$ taming $J|B$, with
$[\omega_B]= c_f|B \in H_{\dR}^2(B)$.
Then $\omega_B$ extends to a closed 2--form $\zeta$ on $X$ representing
$c_f$, with $\nu$
and $TB$ $\zeta$--orthogonal, and $\zeta$ agreeing with the given Hermitian
form on each 1--dimensional $J$--complex subspace of $\nu$.
Given such an extension $\zeta$, suppose that  each $F_y$, $y\in\CP^k$, 
has a neighborhood $W_y$ in $X$ with a closed 2--form $\eta_y$ on $W_y$ taming
$J|\ker df_x$ for all $x\in W_y -B$, agreeing with $\zeta$ on each
$TF_z|B$, $z\in \CP^k$, and
with $[\eta_y-\zeta] =0\in H_{\dR}^2 (W_y,B)$.
Then $(f,J)$ determines an isotopy class $\Omega$ of symplectic forms on $X$
representing $c_f\in H_{\dR}^2 (X)$.
\end{thm}

Each $\ker df_x$ is $J$--complex, so we define $\eta_y$--tameness on it in the 
obvious way. 
The class $[\eta_y-\zeta] \in H_{\dR}^2(W_y,B)$ is defined since
$\eta_y -\zeta$ vanishes on $B$ by hypothesis.
This class vanishes automatically if $[\eta_y] = c_f|W_y$ and the
restriction map $H_{\dR}^1 (W_y )\to H_{\dR}^1 (B)$ is surjective;
however surjectivity always fails when (for example) 
$B$ is a surface of nonzero genus 
and a generic (4--dimensional) fiber has $b_1 <2$.

For our subsequent application to Lefschetz pencils, we will need an explicit
characterization of $\Omega$ and detailed properties of some of its 
representatives. 
The characterization below is complicated by our need to perturb $J$ during 
the proof. 
A simpler version when no perturbation is necessary will be given as 
Addendum~\ref{Jfixed} after the required notation is established.

\begin{adden}\label{unique}
Fix a metric on $X$ and $\ep>0$.
Let $\J_\ep$ be the $C^0$--space of continuous almost-complex structures
$J'$ on $X$ that are $\ep$--close to $J$, agree with $J$ on $TX|B$ and 
outside the $\ep$--neighborhood $U$ of $B$, and make 
each $F_y\cap U$ $J'$--complex.
Fix a regular value $y_0$ of $f$. 
Then $\Omega$ contains a form $\omega$  
taming an element of $\J_\ep$ and 
extending $\omega_B$,  such that 
$J$ is $\omega$--compatible on $\nu$, which is $\omega$--orthogonal to $B$, 
and $\omega|F_{y_0}$ is isotopic to $\eta_{y_0}|F_{y_0}$ 
by an isotopy $\psi_s$ of the pair 
$(F_{y_0},B)$ that is symplectic on $(B,\omega_B)$.
For $\ep $ sufficiently small, any two forms
representing $c_f$ and taming elements  of $\J_\ep$ are isotopic, so
these latter conditions  uniquely determine $\Omega$.
\end{adden}

Theorem~\ref{thm:taming}  was designed for compatibility with 
\cite[Theorem 3.1]{G1}, which was the main tool for putting symplectic 
structures on hyperpencils (and domains of locally holomorphic maps 
\cite{G2}). 
The proof is based on an idea of Thurston \cite{T}. 
We state and prove a version of the theorem which has been slightly 
modified, primarily to correct for the failure of $H^1$--surjectivity 
observed following Theorem~\ref{thm:taming}. 
We will ultimately apply the theorem to a linear system projection 
$f\co X-B\to \CP^k$, working relative to a normal disk bundle $C$ of $B$ 
(intersected with $X-B$). 

\begin{thm} \label{thm:smoothmap}
Let $f\co X \to Y$ be a smooth map between manifolds, 
and let $C$ be a codimension--$0$ submanifold (with boundary) 
that is closed in $X$, with $X-\Int C$ compact. 
Suppose that $\omega_Y$ is a symplectic form on $Y$, and $J$ is a
continuous, $(\omega_Y, f)$--tame almost-complex structure on $X$.
Let $\zeta$ be a closed 2--form on $X$ taming $J$ on $C$. 
Suppose that for each $y \in Y$,
$f^{-1}(y)\cup C$ has a neighborhood $W_y$ in $X$, with 
a closed 2--form $\eta_y$ on $W_y$ agreeing with $\zeta$ on $C$,
such that $[\eta_y-\zeta] = 0 
\in H^2_{\dR}(W_y,C)$ and such that $\eta_y$ tames $J| \ker df_x$
for each $x \in W_y$.
Then there is a closed 2--form
$\eta$ on $X$ agreeing with $\zeta$ on $C$, with $[\eta] = [\zeta] \in
H^2_{\dR}(X)$, and such that for all sufficiently small $t>0$
the form $\omega_t = t\eta + f^*\omega_Y$ on $X$ tames $J$ (and hence is
symplectic).
For preassigned $\hat y_1,\ldots,\hat y_m\in Y$, we can assume $\eta$ agrees 
with $\eta_{\hat y_j}$ near each $f^{-1}(\hat y_j)$.
\end{thm}

\begin{proof} 
For each $y \in Y$, $[\eta_y-\zeta] = 0\in H_{\dR}^2 (W_y,C)$, 
so we can write $\eta_y = \zeta +
d\alpha_y$ for some 1--form $\alpha_y$ on $W_y$ with $\alpha_y|C=0$.
Since each $X-W_y$ is compact, each $y \in Y$ has a neighborhood
disjoint from $f(X-W_y)$.
Thus, we can cover $Y$ by open sets $U_i$,
with each $f^{-1}(U_i)$ contained in some $W_y$,   
and each $\hat y_j$ lying in only one $U_i$. 
Let $\{\rho_i\}$ be a subordinate partition of unity on $Y$.
The corresponding partition of
unity $\{\rho_i \circ f\}$ on $X$ can be used to splice the forms
$\alpha_y$; let $\eta =\zeta + d\sum_i(\rho_i \circ f)\alpha_{y_i}$.
Clearly, $\eta$ is closed with $[\eta] = [\zeta] 
\in H^2_{\dR}(X)$, $\eta = \zeta$ on $C$, 
and $\eta = \eta_{\hat y_j}$ near $f^{-1}(\hat y_j)$, so it suffices to
show that $\omega_t$ tames $J$ ($t > 0$ small).
In preparation, perform the
differentiation to obtain $\eta = \zeta + \sum_i(\rho_i \circ f)
d \alpha_{y_i} + \sum_i (d\rho_i \circ df) \wedge
\alpha_{y_i}$.
The last term vanishes when applied to a  pair of
vectors in $\ker df_x$, so on each $\ker df_x$ we have
$\eta = \zeta + \sum_i(\rho_i \circ f)d\alpha_{y_i} =
\sum_i(\rho_i \circ f)\eta_{y_i}$.
By hypothesis, this is a convex
combination of taming forms, so we conclude that
$J|\ker df_x$ is $\eta$--tame for each $x \in X$.

It remains to show that there is a $t_0 > 0$ for which
$\omega_t(v, Jv) > 0$ for every $t \in (0, t_0)$ and $v$ in the unit
sphere bundle $\Sigma \subset TX$ (for any convenient metric).
But 
$$\omega_t(v, Jv) = t \eta (v, Jv) + f^*\omega_Y(v, Jv).$$
Since $J$ is $(\omega_Y, f)$--tame, the last term is positive for $v \notin
\ker df$ and zero otherwise.
Since $J|\ker df$ is
$\eta$--tame, the continuous function $\eta(v,Jv)$ is positive for all
$v$ in some neighborhood $U$ of $\ker df \cap \Sigma$ in $\Sigma$.
Similarly, for $v \in \Sigma|C$, $\eta(v, Jv) = \zeta (v, Jv) > 0$.
Thus, $\omega_t(v, Jv) > 0$ for all $t > 0$ when $v \in U \cup \Sigma|C$.
On the compact set $\Sigma|(X-\Int C) - U$ containing the
rest of $\Sigma$, $\eta(v,Jv)$ is bounded and the last displayed term is
bounded below by a positive constant, so $\omega_t(v, Jv) > 0$ for $0 <
t < t_0$ sufficiently small, as required.
\end{proof}


\begin{proof}[Proof of Theorem~\ref{thm:taming} and addenda]
We begin by producing the desired symplectic structure near $B$, 
via a local model
generalizing the case $\dim B=0$ from \cite{G1}. 
Assume the fibers of $\pi$ are tangent to  $\nu$.
Let $L_0\to B$ denote the Hermitian line bundle obtained by restricting
$\pi$ to a fixed $F_y$, so $L_0$ and $\nu$ are associated to the same
principal $\rU(1)$--bundle $\pi_P\co P\to B$.
Then $c_1(L_0) = c_f|B = [\omega_B]$ (since a generic section of $L_0$
is obtained by perturbing $B\cup f^{-1}(\CP^{k-1})\subset X$
and intersecting it with $F_y$).
Let $i\beta_0$ on $P$ be a $\rU(1)$--connection form for $L_0$ with Chern form
$\omega_B$, so $-\frac1{2\pi} d\beta_0 = \pi_P^* \omega_B$.
For $r>0$, let $S_r \subset V$ denote the sphere bundle of radius $r$
(for the Hermitian metric).
The map $(\pi,f)\co V-B\to B\times\CP^k$ exhibits each $S_r$ as a
principal $\rU(1)$--bundle.
The corresponding line bundle $L\to B\times \CP^k$ restricts to $L_0$ over
$B$ and to the tautological bundle $L_{\text{taut}}$ over $\CP^k$.
Since $H^2 (B\times\CP^k) \cong H^2(B) \oplus (H^0 (B)\otimes H^2(\CP^k))$
(over $\zed$), we conclude that $L\cong \pi_1^* L_0 \otimes \pi_2^*
L_{\text{taut}}$.
Fix this isomorphism, and let $i\beta$ be the $\rU(1)$--connection form on
$S_r$ induced by $i\beta_0$ on $L_0$ and the tautological connection on
$L_{\text{taut}}$.
Then the Chern form of $i\beta$ is given by $-\frac1{2\pi}d\beta =
\pi^* \omega_B - f^* \ostd$ (pushed down to $B\times \CP^k$).
Define a 2--form $\omega_V$ on $V-B$ by
$$\omega_V = (1-r^2) \pi^* \omega_B + r^2 f^* \ostd
+ \frac1{2\pi} d(r^2)\wedge \beta.$$
An easy calculation shows that $d\omega_V =0$, and it is routine to verify 
\cite{G1} 
that $\omega_V$ restricts to
the given Hermitian form on each fiber of $\pi$ (up to a constant factor 
of $\pi$, arising from our choice of normalization of $\ostd$, which can 
be eliminated by a constant rescaling of $r$). 
Let $H$ be the smooth distribution on $V$ consisting of $TB$ on $B$
together with its $\beta$--horizontal lifts to each $S_r$.
Clearly, $H$ is tangent to each $S_r$ and $F_y$, so it is
$\omega_V$--orthogonal to the fibers of $\pi$.
Since $\omega_V |H = (1-r^2)\pi^* \omega_B$ extends smoothly over $B$,
as does $\omega_V$ on the $\pi$--fibers, $\omega_V$ extends smoothly to
all of $V$, with $\omega_V|B=\omega_B$.
If $J_V$ denotes the almost-complex structure on $V$ obtained by lifting
$J|B$ to $H$ and summing with the complex bundle structure on the fibers of
$\pi$, then $J_V$ is $\omega_V$--tame for $r<1$.
(Check this separately on the $\pi$--fibers and their $\omega_V$--orthogonal
complements $H$.)
Note that $J_V =J$ on $TX|B$ (Lemma~\ref{nu}).

We can now state the remaining addendum:

\begin{adden}\label{Jfixed}
If $J$ agrees with $J_V$ near $B$ for some choice of 
$\pi\co V\to B$ and $\beta_0$ as above, then
$\Omega$ has the simpler characterization that it contains forms $\omega$
taming $J$ with $[\omega]=c_f$.
In fact, there is a $J$--taming form $\omega\in\Omega$ 
satisfying all the conclusions of Addendum~\ref{unique}
with $\nu$ induced by $\pi$, and such that the given forms 
$\psi_s^*\eta_{y_0}$ on $F_{y_0}$ between $\eta_{y_0}$ and $\omega$ 
all tame $J$.
\end{adden}

To construct the required form $\zeta$,
choose a form $\zeta_0$ representing $c_f \in H_{\dR}^2 (X)$.
Then $[\omega_V -\zeta_0] = 0\in H_{\dR}^2 (V)$, so there is a 1--form
$\alpha$ on $V$ with $d\alpha = \omega_V -\zeta_0$.
Let $\zeta = \zeta_0 + d(\rho\alpha)$, where $\rho\co X\to \real$ has
support in $V$ and $\rho=1$ near $B$.
Then $\zeta = \omega_V$ near $B$, so $\zeta$ satisfies the required
conditions for the theorem.
If $\zeta_0$ was already the hypothesized extension of $\omega_B$, satisfying these conditions and
suitably compatible with forms $\eta_y$, then  $\zeta_0 = \omega_V=\zeta$
on each $TF_z|B \cong TB \oplus L_0$, 
so $\zeta$  still agrees with each $\eta_y$ as required along $B$.
We also could have arranged $\alpha|B=0$ since $H_{\dR}^2(V,B)=0$, 
so that we still have $[\eta_y-\zeta] =0 \in H_{\dR}^2 (W_y,B)$.
Thus, we can assume the given $\zeta$ agrees with $\omega_V$ near $B$.

Since we must perturb $J$ near $B$, 
we verify that for sufficiently small $\ep$, 
every $J'\in \J_\ep$ as in Addendum~\ref{unique}  
is $(\omega_{\std},f)$--tame on $X-B$.
Choose $\ep$ so that the
$\ep$--neighborhood $U$ of $B$ in $X$ (in the given metric) has closure
in $V$, and let $\Sigma \subset TX$ be the compact subset consisting
of unit vectors over $\cl(U)$ that are $\omega_V$--orthogonal to fibers
$F_y$.
For $J' \in \J_\ep$, each $\ker df_x = T_x F_{f(x)}$ over $U-B$ is
$J'$--complex, so it suffices to show that $f^*\omega_{\text{std}}(v,
J'v) > 0$ for $v \in \Sigma\cap T(U-B)$.
We replace $f^*\omega_{\text{std}}$ by $\omega_V$, since these agree on
such vectors $v$ 
(which are tangent to the $\pi$--fibers and $S_r)$ 
up to the scale factor $r^2 > 0$.
But $\omega_V(v,Jv)>0$ for $v \in \Sigma$ (since $J$ equals $J_V$ on
$TX|B$ and $J$ is $(\omega_{\text{std}}, f)$--tame elsewhere), so the
corresponding inequality holds for all $J' \in \J_\ep$ for $\ep$
sufficiently small, by compactness of $\Sigma$ and openness of the taming 
condition. 

We must also modify the pairs $(W_y,\eta_y)$
so that for all sufficiently small $\ep$, 
every $J'\in \J_\ep$ is 
$\eta_y$--tame on $\ker df_x$ for each $y\in\CP^k$ and $x\in W_y-B$.
Shrink each $W_y$ so that $\eta_y$ is defined on  $\cl(W_y)$.
Each $W_y$ contains $f^{-1}(U_y)$ for some neighborhood
$U_y$ of $y$  (cf proof of Theorem~\ref{thm:smoothmap}).
After passing to a finite subcover of $\{U_y\}$, we can assume $\{W_y\}$ is
finite, so the pairs $(W_y,\eta_y)$ for all $y\in \CP^k$ are taken from
a finite set, and $\eta_{y_0}|F_{y_0}$ is preserved.
Now for each $\eta_y$, $\eta_y(v, Jv) > 0$ on the compact space of
unit tangent vectors to fibers $F_y$ in $\cl (W_y\cap V)$.
(Note that on $TX|B$, $\zeta$ tames $J$.)
Thus,  each $J'\in\J_\ep$ has
the required $\eta_y$--taming for $\ep$ sufficiently small.

Next we splice our local model $\omega_V$ and $J_V$ into each $\eta_y$ and $J$.
For $y\in \CP^k$, $\eta_y$ equals $\zeta$ on $TF_y|B$, so it tames
$J$ there and hence is symplectic on $F_y$ near $B$.
Thus, Weinstein's symplectic tubular neighborhood theorem 
\cite{W} on $F_y$ produces an isotopy of $F_y$ fixing $B$ (pointwise) 
and supported in a preassigned neighborhood of $B$, sending
$\eta_y|F_y$ to a form $\eta'_y$ 
agreeing with $\zeta= \omega_V$ near $B $ on $F_y$.
To extend $\eta'_y$ to a neighborhood of $F_y$ in $X$, 
first extend it as $\omega_V$ near $B$ and as $\eta_y$ farther away, 
leaving a gap in between (inside $V$). 
Let $r\co W_y\to W_y$ be a smooth map agreeing with $\id_{W_y}$ away from the 
gap and on $F_y$, collapsing $W_y$ onto $F_y$ near the gap. 
Then $r^* \eta'_y$ is a closed form near $F_y$  extending $\eta'_y$ 
(cf \cite{G1}). 
Now recall that the  vector field generating Weinstein's isotopy 
vanishes to second order on $B$. 
(It is symplectically dual to the 1--form $-\int_0^1\pi_t^* (X_t\contract
(\eta_y-\zeta))\,dt$, where $\pi_t$ is fiberwise multiplication by $t$, 
and the radial vector field $X_t = \frac{d}{dt}\pi_t$ vanishes to first 
order on $B$, as does $\eta_y -\zeta$.) 
Thus we can
assume our isotopy is arbitrarily $C^1$--small (by working in a
sufficiently small neighborhood of $B$),
so we can replace $\eta_y$
on $W_y$ by $r^* \eta'_y$ on a sufficiently small neighborhood of $F_y$
without disturbing our original hypotheses.
In particular, we can assume $J$ is still $\eta_y$--tame on each
$\ker df_x$ (or similarly for all $J'$ in a preassigned compact
subset of $\J_\ep$ with $\ep$ as in the previous paragraph).
Since we have shrunk the sets $W_y$, the set $\{W_y\}$ may again be
infinite, but we can reduce to a finite subcollection as before.
Then there is a single neighborhood $W$ of $B$ in $X$, contained in
$\bigcap W_y$, on which each $\eta_y$ agrees with $\omega_V$ and $\zeta$.
Since $\eta_{y_0}|F_{y_0}$ has only been changed by a $C^1$--small isotopy 
fixing $B$, its use in the addenda is unaffected. 

To complete splicing the local model, we perturb $J$ to $J'$ agreeing 
with $J_V$ near $B$. 
Under the hypothesis of Addendum~\ref{Jfixed}, we simply set $J'=J$. 
Otherwise, we invoke \cite[Corollary 4.2]{G1}, which was adapted from 
\cite[page 100]{ABKLR}. 

\begin{lem}\label{retract}
\cite{G1}\qua 
For any finite dimensional, real vector space $V$, there is a canonical 
retraction $j(A) = A(-A^2)^{-1/2}$ from the open subset of operators in 
$\Aut (V)$ without real eigenvalues to the set of 
linear complex structures on $V$.
For any linear $T\co V\to W$ with $TA=BT$, we have $Tj(A) = j(B)T$
(when both sides are defined). 
\end{lem}

Since $J=J_V$ on $TX|B$, $J_t = j((1-t)J+tJ_V)$ is well-defined for 
$0\le t\le1$ near $B$, and each $F_y$ is $J_t$--complex there (as seen by 
letting $T$ be inclusion $T_xF_y\to T_xX$). 
For any $\ep>0$, we can thus define $J'\in \J_\ep$ to be $J_\rho$, for 
$\rho\co X\to I$ supported sufficiently close to $B$ and with $\rho\equiv 1$ 
near $B$, extended as $J$ away from $\supp \rho$.
Then for $\ep$ sufficiently small, the preceding three paragraphs show that 
$(f,J')$ is a linear $k$--system satisfying the hypotheses of 
Theorem~\ref{thm:taming} with $J',\zeta$ and each $\eta_y$ agreeing with 
the standard model on a suitably reduced $W$.

We now construct a symplectic form $\omega$ on $X$ as in \cite{G1}. 
First we apply Theorem~\ref{thm:smoothmap} to $f\co X- B\to \CP^k$
and $J'$, with $C\subset W$ a normal disk bundle to $B$ (intersected
with $X-B$).
Note that $[\eta_y -\zeta] \in H_{\dR}^2 (W_y - B,C) \cong H_{\dR}^2 (W_y,B)$ 
vanishes as required. 
We obtain a closed 2--form 
$\eta$ on $X-B$ agreeing with $\omega_V$ on $C$ (hence extending 
over $X$), with $[\eta] = c_f\in H_{\dR}^2 (X)$ and $\eta = \eta_{y_0}$ 
on $F_{y_0}$, such that $\omega_t = t\eta + f^* \ostd$ tames $J'$ on 
$X-B$ for $t>0$ chosen sufficiently small. 
On $C$, the symplectic form $\omega_t$ is given by 
$$\omega_t (r) = t(1-r^2)\pi^* \omega_B + (1+tr^2) f^* \ostd + 
\frac{t}{2\pi} d(r^2) \wedge\beta.$$
Unfortunately, this is singular at $B$. 
(Compare the middle term with that of $\omega_V\ne 0$.) 
However, we can desingularize by a dilation in the manner of \cite{G1}:
The radial change of variables $R^2 =
\frac{1+tr^2}{1+t}$  shows that $\omega_V(R) = \frac{1}{1+t}
\omega_t(r)$, so there is a radial symplectic embedding
$\varphi\co(C,\frac{1}{1+t}\omega_t) \to (V,\omega_V)$
onto a collar surrounding the bundle $R^2 \le \frac{1}{1+t}$.
Let $\varphi_0\co V \to V$ be a radially symmetric diffeomorphism 
covering $\id_B$ and agreeing with $\varphi$ near $\partial C$. 
Let $\omega$ be $\varphi_0^*\omega_V$ on $C\cup B$ and 
$\frac{1}{1+t} \omega_t$ elsewhere.
These pieces fit together to define a symplectic form on $X$, 
since $\varphi$ is a symplectic embedding. 
(This construction is equivalent to blowing up $B$,
applying Theorem~\ref{thm:smoothmap} with $C = \emptyset$ to the resulting
singular fibration, and then blowing back down, but it avoids
technical difficulties associated with taming on the blown up base locus.)

The form $\omega$ satisfies the properties required by 
Theorem~\ref{thm:taming} and its addenda: 
To compute the cohomology class $[\omega] \in
H^2_{\dR}(X)$, it suffices to work outside $C$.  
Then $[\omega] = \frac{1}{1+t}[\omega_t] = \frac{1}{1+t}(t c_f +
f^*[\omega_{\std}]) = c_f$ as required, since $[\omega_{\std}] = h
\in  H^2_{\dR}(\CP^k)$. 
For Addendum~\ref{unique}, note that $\omega$ obviously extends $\omega_B$
and is compatible with $J$ on $\nu$, which is $\omega$--orthogonal to $B$. 
Outside $C$, we already know that $\omega = \frac1{1+t} \omega_t$ tames 
$J'\in \J_\ep$, so taming need only be checked for $J'= J_V$ on $C\cup B$ with 
$\omega = \varphi_0^* \omega_V$, and this is easy on $TX|B=TB\oplus \nu$. 
For $C$, consider the $\omega_V$--orthogonal, $J_V$--complex splitting 
$T(V-B) = P\oplus P^\bot
\oplus H$, where $P$ and $P^\bot$ are tangent and normal, respectively,
to the complex lines through $B$ in the bundle structure $\pi$.
The radial map $\varphi_0$ preserves the splitting but scales each summand 
by a different positive function. 
(Although the fibers of $P$ are scaled differently along their two axes, 
$\varphi_0^*$ only rescales $\omega_V|P$ since it is an area form.)
Now $J_V$ is $\omega$--tame on $C$ since it is $\omega_V$--tame on each 
summand. 
To verify that $\omega|F_{y_0}$ is pairwise isotopic to $\eta_{y_0}|F_{y_0}$,
recall that $\eta|F_{y_0}= \eta_{y_0}|F_{y_0}$, so
$\omega |F_{y_0} = \frac{t}{1+t} \eta_{y_0}|F_{y_0}$ outside $C$. 
When $t\to\infty$ we have $R\to r$ and $\varphi_0 \to \id_V$, 
so $\omega\to\eta$. 
Note that $\eta$ and $\omega$ (for all $t>0$) are symplectic on $F_{y_0}$, 
although not necessarily on $X$ (unless $t$ is small). 
The required isotopy now follows from Moser's method \cite{M} applied 
pairwise to $(F_{y_0},B)$: 
Starting from $\omega$ as constructed above with $t$ sufficiently small, 
let $\widetilde\omega_s$, $s= \frac1t \in [0,a]$, be the corresponding family 
of cohomologous symplectic forms on $F_{y_0}$ obtained by letting 
$t\to\infty$ (so $\widetilde\omega_a = \omega|F_{y_0}$ and 
$\widetilde\omega_0 = \eta_{y_0}|F_{y_0}$). 
Moser gives a family $\alpha_s$ of 1--forms on $F_{y_0}$ 
with $d\alpha_s = \frac{d}{ds}\widetilde{\omega}_s$,
then flows by the vector field $Y_s$ for which $\widetilde{\omega}_s(Y_s,\cdot)=
- -\alpha_s$ to obtain an isotopy with $\psi_s^*\eta_{y_0} =\widetilde\omega_s$.
If we first subtract $dg_s$ from $\alpha_s$, where $g_s\co F_{y_0}\to\real$ is
obtained by pushing $\alpha_s\co TF_{y_0}\to \real$ from 
$TB^{\bot_{\widetilde\omega_s}}$ down to a tubular
neighborhood of $B$ and tapering to $0$ away from $B$, then we can assume
$\alpha_s |TB^{\bot_{\widetilde\omega_s}} =0$. 
Thus $Y_s$ is $\widetilde\omega_s$--orthogonal to 
$TB^{\bot_{\widetilde\omega_s}}$, so  $Y_s$ is tangent to $B$, and its flow 
$\psi_s$ preserves $B$ as
required, completing verification of the conditions of Addendum~\ref{unique}. 
(The isotopy restricts to symplectomorphisms on $B$ since each 
$\widetilde\omega_s|B=\omega_B$.) 
Addendum~\ref{Jfixed} now follows immediately from the observation that 
the forms $\widetilde\omega_s= \psi_s^* \eta_{y_0}$ on $F_{y_0}$ all tame 
$J = J'$ in this case.
(For the characterization of $\Omega$, note that any two cohomologous forms
taming a fixed $J$ are isotopic by convexity of the taming condition and 
Moser's Theorem.) 

To complete the proof of Theorem~\ref{thm:taming} and Addendum~\ref{unique}, 
we show that for sufficiently small $\delta$, any two forms $\omega_u$, 
$u=0,1$, taming structures $J_u\in \J_\delta$
and representing $c_f \in H_{\dR}^2 (X)$, are isotopic, implying that
$\Omega$ is canonically defined. 
(Note that for $0<\delta <\ep$ we have $J\in \J_\delta \subset \J_\ep$, 
so $\Omega$ is then independent of sufficiently small $\ep>0$ and agrees 
with its usage in Addendum~\ref{Jfixed}. 
Metric independence follows, since for metrics $g,g'$ on $X$ and $\ep>0$ 
there is a $\delta >0$ with $\J_\delta (g')\subset \J_\ep(g)$.)
Let $J_u = j((1-u)J_0 + uJ_1)$, $0\le u\le 1$.
For $\delta$ sufficiently small, this is a well-defined path from
$J_0$ to $J_1$, and each $J_u$ satisfies the defining conditions for 
$\J_\delta$ except possibly for $\delta$--closeness to $J$. 
For $\delta$ sufficiently small, there is a compact subset $K$ of the 
bundle $\Aut (TX)\to X$ lying in the domain of $j$, containing a 
$\delta$--neighborhood of the image of the section $J$. 
By uniform continuity of $j|K$, we can
choose $\delta \in (0,\ep)$ such that $J_u$ must be a path in $\J_\ep$,
with $\ep$ small enough to satisfy all of the previous requirements.
Now for fixed $J_0,J_1\in\J_\delta$, 
we can assume the forms $\eta_y$ were constructed as 
above to agree with $\omega_V$ on $W$ and tame each 
$J_u|\ker df_x$, $0\le u\le1$.
Perturb the entire family as before to $J'_u \in \J_\ep$, $0\le u\le1$, with
each $J'_u$ agreeing with $J_V$ on a fixed $W$.
For a small enough perturbation, $J'_u$ will be $\omega_u$--tame, $u=0,1$.
For $0<u<1$, the previous argument produces symplectic forms $\omega_u$
taming $J'_u$.
The family $\omega_u$, $0\le u\le1$, need not be continuous. 
However, each $\omega_u$ tames $J'_v$ for $v$ in a neighborhood of $u$, 
so splicing by a partition of unity on the interval $I$ produces (by convexity
of taming) a smooth family $\omega'_u$ taming $J'_u$, $0\le u\le1$, 
with $\omega'_u = \omega_u$ for $u=0,1$. 
Applying Moser's Theorem to this family of cohomologous symplectic 
forms gives the required isotopy. 
\end{proof}
\newpage

\section{Lefschetz pencils}	

We now return to the investigation of Lefschetz pencils 
(Definition~\ref{de:pencil}) and complete the discussion of their 
classification theory (Proposition~\ref{prop:fixedA}, cf Principle~\ref{B}). 
We then apply the results of the previous section on linear 1--systems, 
to show that a similar topological classification theory applies in the 
symplectic setting (Theorem~\ref{thm:omega}, cf Principle~\ref{C}). 

To analyze the topology of a Lefschetz critical point (eg \cite{L}), 
recall the local model ${f\co\cee^n\to\cee}$, $f(z) = \sum_{i=1}^n z_i^2$, 
given in Definition~\ref{de:pencil}(2). 
To see that a 
regular neighborhood of the singular fiber is obtained from that of a
regular fiber by adding an $n$--handle, note that the 
core of the $n$--handle appears in the local model 
as the $\ep$--disk $D_\ep$ in $\real^n\subset \cee^n$.
Thus, the handle is attached to the fiber $F_{\ep^2}$ along an embedding
$S^{n-1}\hookrightarrow F_{\ep^2} -B$ whose normal bundle
$\nu S^{n-1}=-iTS^{n-1}$ in the complex bundle $TF$ 
is identified with $T^*S^{n-1}$ (by contraction with
$\omega_{\cee^n}$).
We will call such an embedding, together with its isomorphism
$\nu S^{n-1}\cong T^* S^{n-1}$, a {\em vanishing cycle}.
Regular fibers intersect the local model in manifolds diffeomorphic
to $T^* S^{n-1}$, and the singular fiber is obtained by collapsing
the 0--section (vanishing cycle) to a point.
(The latter assertion can be seen explicitly by writing the real and
imaginary parts of the equation $\sum z_i^2 =0$ as $\|x\| = \|y\|$,
$x\cdot y=0$.)
The monodromy around the singular fiber is obtained from the
geodesic flow on $T^*S^{n-1}\cong TS^{n-1}$, renormalized to be
$2\pi$--periodic near the 0--section (on which the flow is undefined), 
and tapered to have compact support \cite{Ar,S2}.
At time $\pi$, the resulting diffeomorphism extends over the 0--section
as the antipodal map, defining the monodromy, which is called a (positive)
Dehn twist.
(To verify this description, note that multiplication by $e^{i\theta}$
acts as the $2\pi$--periodic geodesic flow on the singular fiber, and
makes $f$ equivariant with respect to $e^{2i\theta}$ on the base.
Thus the sphere $\partial D_\ep$ is transported around the singular
fiber by $e^{i\theta}$,
returning to its original position when $\theta = \pi$, with antipodal
monodromy.
Away from $\partial D_\ep$, the
monodromy is obtained via projection to the singular fiber, where it
can be tapered from the geodesic flow near 0 to the identity outside a
compact set by an isotopy.)
Given arcs $A=\bigcup A_j$ in $\CP^1$ as in the introduction, 
connecting each critical value of a Lefschetz pencil
to a fixed regular value,
say [1:0], we may  interpret all vanishing cycles and
monodromies as occurring on the single fiber $F_{[1:0]}$.
The disk $D_\ep$ at each critical point extends to a disk $D_j$ with
$f(D_j)=A_j$ and $\partial D_j$ the vanishing cycle in $F_{[1:0]}$.
Following Lefschetz, we will call such a disk a {\em thimble},
but we also require that each $f|D_j\co D_j\to A_j$ has a nondegenerate, 
unique critical point, and that there is a local trivialization of $f$
near $F_{[1:0]}$ in which each $D_j$ is horizontal.

All of the above  structure on the local model of a critical point is
compatible with suitable symplectic forms.
(See \cite{S2}.)
For the standard K\"ahler form $\omega_{\cee^n}$,
the sphere $\partial D_\ep$ in $\cee^n$ 
is Lagrangian in the symplectic submanifold
$F_{\ep^2}$, so by Weinstein's theorem \cite{W} it has a neighborhood
symplectomorphic to a neighborhood of the 0--section in $T^*S^{n-1}$.
This allows Dehn twists to be defined symplectically by a Hamiltonian
flow in $T^* S^{n-1}-(0\text{--section})$ \cite{Ar,S2},
determining the monodromy
around the singular fiber up to symplectic (Hamiltonian) isotopy.
The Lagrangian embedding $S^{n-1} \hookrightarrow F_{\ep^2}-B$
determines a vanishing cycle,
and will be called the {\em Lagrangian vanishing cycle}
for the critical point.
If $\omega$ is an arbitrary K\"{a}hler form near 0 on $\complex^n$,
any small arc $A_j$ from $0\in\complex$ (such as $[0,\ep^2]$ above) still
determines a smooth Lagrangian thimble and vanishing cycle, by a trick of
Donaldson \cite[Lemma 1.13]{S2}.
The disk consists of the trajectories under symplectic parallel transport
(ie the flow over $A_j$ $\omega$--normal to the fibers)
that limit to the critical point.
If $\omega$ is only given to be compatible with $i$ at 0, this
structure still exists.
(In fact, $\omega$ agrees at 0 with some K\"{a}hler form; after
rescaling the coordinates, we may assume the two forms are arbitrarily
close, as are the resulting disks and vanishing cycles.
The case of an arbitrary taming $\omega$ is less clear.)
For a given Lefschetz pencil $f\co X-B \to \CP^1$, arcs $A \subset \CP^1$,
and symplectic form $\omega$ on $X$ that is symplectic on each $F_y-K$ 
(where $K\subset X-B$ is the critical set as in 
Definition~\ref{de:pencil}(3)), any
such disk at $x \in K$ is uniquely determined and uniquely extends to
a Lagrangian thimble, by symplectic parallel transport.

For a closed, oriented manifold pair $B\subset F$ of dimensions $2n-4$
and $2n-2$, respectively,
let $\D = \D(F,B)$ denote the group of orientation-preserving
self-diffeomorphisms of $F$ fixing (pointwise) $B$ and $TF|B$.
If $\omega_F$ is a symplectic form on $F$ whose restriction to $B$
is symplectic, let $\D_{\omega_F} = \D_{\omega_F} (F,B) \subset \D$ be the
subgroup of symplectomorphisms of $F$ fixing $B$ and $TF|B$.
Let $\delta \in \pi_0 (\D)$ be the element obtained by a $2\pi$
counterclockwise rotation of the normal fibers of $B$, extended
in the obvious way (by tapering to $\id_F$) to a diffeomorphism of $F$.
In the symplectic case, $\delta$ canonically pulls back to
$\delta_{\omega_F} \in \pi_0(\D_{\omega_F})$, determined by the Hamiltonian
flow of a suitable radial function on a tubular neighborhood of $B$.
Any vanishing cycle in $F-B$ determines a Dehn twist in $\pi_0(\D)$
as described above.
In the symplectic case, a Lagrangian embedding $S^{n-1}\hookrightarrow F-B$
determines a symplectic Dehn twist in $\pi_0(\D_{\omega_F})$, whose image
in $\pi_0(\D)$ is generated by the corresponding vanishing cycle.
If $\omega'$ is obtained from $\omega_F$ 
by a pairwise diffeomorphism of $(F,B)$,
there is an induced isomorphism $\D_{\omega'}\cong\D_{\omega_F}$
sending $\delta_{\omega'}$ to $\delta_{\omega_F}$
and inducing an obvious correspondence of symplectic Dehn twists.

\begin{prop}\label{prop:fixedA}
Let $w= (t_1,\ldots,t_m)$ be a word in positive Dehn twists
$t_j\in \pi_0 (\D)$, whose product $\prod_{j=1}^m t_j$ equals $\delta$.
Then there is a manifold $X$ with a Lefschetz pencil $f$ whose fiber
over $\text{\rm [1:0]} \in \CP^1$ is $F\subset X$, whose base locus is $B$, and
whose monodromy around the singular fibers (with respect to fixed
arcs $A\subset\CP^1$) is given by $w$.
For a fixed choice of $A$ and vanishing cycles determining the Dehn twists,
such Lefschetz pencils are classified by $\pi_1(\D)$.
\end{prop}

\begin{proof}
Choose $0<\theta_1 <\cdots <\theta_m<2\pi$.
For each $j$, attach an $n$--handle to $D^2\times F$ using the given
vanishing cycle for $t_j$ in $\{e^{i\theta_j}\}\times F$.
We obtain a singular fibration over $D^2$, with $m$ singularities as
described above, and monodromy given by $w$.
Since $\prod_{j=1}^m t_j =\delta$ is isotopic to $\id_F$ fixing  $B$ (but
rotating its normal bundle clockwise), the fibration over $\partial D^2$ can be
identified with $\partial D^2 \times (F,B)$, and the freedom to choose
this identification (without losing control of $TF|B$) is given by $\pi_1(\D)$.
For any such identification, we can glue on a copy of $D^2\times F$ to obtain
a Lefschetz fibration $\widetilde f\co\widetilde X \to \CP^1$.
This $\widetilde{X}$ contains a canonical copy of $\CP^1 \times B$,
on which $\widetilde f$ restricts to the obvious projection.
The twist defining $\delta$ forces the normal bundle of $\CP^1\times B$
to restrict to the tautological bundle on each $\CP^1\times \{b\}$, so
we can blow down the submanifold to obtain the required Lefschetz pencil.
\end{proof}

To completely determine the above correspondence between Lefschetz pencils 
and $\pi_1(\D)$, we must make a choice determining which pencil corresponds 
to $0\in \pi_1(\D)$. 
For a fixed Lefschetz pencil, arcs $A$ and vanishing cycles, assume the disk 
$D\subset \CP^1$ containing $A$ is embedded so that $1\in \partial D$ maps 
to the central vertex $[1\colon\! 0]$ of $A$.
The monodromy around $\partial D$ is given to us as a product of Dehn twists,
each of which is well-defined up to isotopies supported near its 
vanishing cycle. 
We choose an arc $\gamma$ in $\D$ from this product to the rotation
determining $\delta$. 
(We are given that such arcs exist.) 
Since the rotation untwists to $\id_F$ by a canonical isotopy of $F$ fixing 
$B$, we have now fixed an identification of the fibration over $\partial D$ 
with $\partial D\times (F,B)$, determining the correspondence with 
$\pi_1(\D)$. 
Note that the freedom to change $\gamma$ is essentially $\pi_1(\D)$, so 
unless we fix $\gamma$, the correspondence is only determined up to 
translations in $\pi_1(\D)$. 
In the symplectic setting, we choose $\gamma$ in $\D_{\omega_F}$ similarly, 
to fix a correspondence with $\pi_1(\D_{\omega_F})$.
In this case, passing back to the smooth setting results in a 
correspondence between pencils and $\pi_1(\D)$ that changes with our choice 
of $\gamma$ in $\D_{\omega_F}$ only through translation by elements of 
$\Im i_*$, where $i_*\co\pi_1 (\D_{\omega_F}) \to \pi_1 (\D)$ is induced 
by inclusion. 
In particular, $\omega_F$ picks out a subcollection of pencils 
corresponding to $\Im i_*$ that is independent of our choice of $\gamma$ 
in $\D_{\omega_F}$. 
We will see that these are precisely the pencils admitting symplectic 
structures suitably compatible with $\omega_F$. 
For our symplectic classification, we wish to allow some flexibility 
in the form over the model fiber $F= F_{[1:0]}$, so we only require it 
to be suitably isotopic to $\omega_F$. 
However, the subtlety in specifying the correspondence with 
$\pi_1 (\D_{\omega_F})$ forces us to keep track of a preassigned isotopy on $F$.
Thus, we classify pairs consisting of a suitable symplectic form $\omega$ 
on $X$ and a suitable isotopy from $\omega| F_{[1:0]}$ to $\omega_F$, 
up to deformations of such pairs. 

To state the theorem we need one further fact. 
It is natural to study symplectic forms on $X$ that are symplectic on the 
fibers of $f$, but this condition makes no sense on the critical set $K$. 
For that we show that $f$ determines a complex structure $J^*$ on $TX|K$, 
and require our forms to be compatible with $J^*$. 
We also require similar compatibility normal to $B$.

\begin{lem}\label{Jstar}
A Lefschetz pencil canonically determines a complex structure $J^*$ on 
$TX|K$ (for $n\ne1$) and on any subbundle $\nu$ of $TX|B$ complementary 
to $TB$. 
$J^*$ is obtained by restricting any
local $(\omega_{\text{std}}, f)$--tame almost-complex structure $J$ 
defined near a point of $K$ or $B$, provided $\nu$ is $J$--complex in the 
latter case. 
\end{lem}

\begin{proof}
First check that in a standard chart at $x \in K$, each hyperplane through
$0$ is a limit of tangent spaces to regular fibers,
so it is $J$--complex for any $(\omega_{\text{std}},f)$--tame local $J$.
Any 1--dimensional complex subspace at $x$ 
is an intersection of such hyperplanes,
so it is also $J$--complex for any such $J$.
But $J_x$ is uniquely determined by its complex lines for $n\ne 1$ 
(\cite[Lemma 4.4(a)]{G1}, cf also Lemma~2.2). 
For $x\in B$, we obtain a suitable $J$ from the complex bundle structure 
$\pi$ of Definition~\ref{de:pencil}(1), 
by perturbing the latter to have fibers tangent to $\nu$ as 
preceding Lemma~\ref{nu}.
That lemma (which only requires $J$ locally) then gives uniqueness on $\nu$.
\end{proof}

\begin{thm}\label{thm:omega} 
Let $\omega_F$ be a symplectic form on $(F,B)$ as preceding 
Proposition~\ref{prop:fixedA}, with 
$[\omega_F]\in H_{\dR}^2 (F)$ Poincar\'e dual to $B$.
Let $S_1, \ldots, S_m$ be Lagrangian embeddings $S^{n-1}
\hookrightarrow F-B$ determining a
word $w = (t_1,\ldots,t_m)$ in positive
symplectic Dehn twists with $\prod_{j=1}^m t_j=
\delta_{\omega_F} \in \pi_0(\D_{\omega_F})$. 
If $n=2$, assume each component of each $F-S_j$ intersects $B$.
Then the corresponding symplectic
Lefschetz pencils are classified by $\pi_1(\D_{\omega_F})$.
More precisely, a Lefschetz pencil $f\co X-B \to \CP^1$ obtained from
$S_1,\ldots,S_m$ as in Proposition \ref{prop:fixedA}, with a fixed choice 
of thimbles $D_j$ bounded by $S_j$ and covering the given arcs $A$,
corresponds to an element of $\Im i_*$
if and only if $X$ admits a symplectic structure $\omega$ that
\begin{enumerate}
\item[\rm(1)] on $(F_{[1:0]},B)$ comes with a pairwise isotopy  to $\omega_F$,
defining a deformation of forms that is fixed on $B$ and $S_1,\ldots,S_m$,
\item[\rm(2)] is symplectic on each $F_y-K$ and (for $n\ne 2$) Lagrangian on each $D_j$,
\item[\rm(3)] is compatible with 
$J^*$  on $TX|K$ (for $n\geq 2$) and on the $\omega$--normal bundle $\nu$ of
$B$,  and
\item[\rm(4)] satisfies $[\omega] = c_f\in H_{\dR}^2(X)$.
\end{enumerate}
For fixed $f$ and $D_1,\ldots,D_m$,
such forms are classified up to deformation through such forms by 
$\pi_2 (\D/\D_{\omega_F})$, and classified by $\ker i_*$ 
if symplectomorphisms preserving $f$ and fixing $f^{-1}(A)$ are also allowed.
\end{thm}

Note that in this classification, a single $\omega$ with different
isotopies as in (1) could represent distinct equivalence classes.
Since a deformation with $[\omega]$ fixed determines an isotopy, the
isotopy classes of forms on $X$ as above 
(for fixed $f$) are classified by the 
quotient of $\pi_2(\D/\D_{\omega_F})$ by
some equivalence relation.
In our case, $[\omega]= c_f$ is Poincar\'e dual to the fiber class 
$[F_{[1:0]}]$, the same condition that arises in Donaldson's construction 
of Lefschetz pencils on symplectic manifolds \cite{D}. 

\begin{proof}
First we assume $n\ge3$ and prepare to apply Theorem~\ref{thm:taming} by
constructing a smooth family $\sigma$ of symplectic structures on the
fibers of a fixed $f$.
As in that proof, the cohomology class of $\omega_B =\omega_F|B$
equals the normal Chern class of $B$ in $F$, so we can define the model
symplectic form $\omega_V$ and almost-complex structure $J_V$ on $V
\subset X$
as before (starting from any fixed choice of $\pi\co V\to B$ as in
Definition~\ref{de:pencil}(1) and any $\omega_B$--tame $J_B$ on $B$).
By Weinstein's theorem, we can assume $\omega_V$ agrees with
$\omega_F$ near $B$
on $F= F_{[1:0]}$.
At each critical point $x_j$, choose a standard chart 
for $f$ (necessarily inducing $J^*$ on $T_{x_j}X$). 
Then $D_j$ is not tangent to any complex curve at $x_j$ (since $f|D_j$ 
is nondegenerate and $f$ is constant or locally surjective on any 
complex curve). 
Thus there is a complex isomorphism $(T_{x_j} X,T_{x_j}D_j)\cong 
(\cee^n ,\real^n)$, and $\omega_{\cee^n}$ pushes down to a symplectic 
form $\omega_j$ near $x_j$, compatible with $J^*$ at $x_j$, and 
Lagrangian on a disk $\Delta_j \subset D_j$ containing $x_j$. 
Since $D_j$ is a thimble for $f$, we can
identify the fibers over $\Int A_j$ with $F$ by an isotopy  in $X$
fixing $B$ and preserving $D_j$, so that $S_j$
matches with $\partial \Delta_j$ (with the correct normal correspondence)
and their tubular neighborhoods in the fiber correspond symplectically
(by Weinstein) relative to $\omega_F$ and $\omega_j$, respectively.
We can assume (by $\rU(2)$--invariance of $\omega_V$) that the isotopy is
$\omega_V$--symplectic near $B$, and that it agrees near $K$ with symplectic
parallel transport in the local model.
Thus it maps the tubular neighborhood of $\partial \Delta_j$ in
its fiber to a neighborhood of the singularity in
the singular fiber, by a map that is a symplectomorphism except on
$\partial \Delta_j$,
which collapses to the singular point (cf \cite{S2}).
Pulling $\omega_F$ back by the isotopy now gives a family $\sigma$ of
symplectic structures on the fibers over $A$,
agreeing with the local models $\omega_V$ and $\omega_j$  near $B$ and $K$.
Extend $\sigma$ over a disk $D \subset \CP^1$ whose interior contains
$A - \text{\rm [1:0]}$.
Since $\prod_{j=1}^m t_j = \delta_{\omega_F}$
is symplectically isotopic to $\id_F$ fixing $B$ (rotating its
normal bundle), we can fix a path $\gamma$ in $\D_{\omega_F}$ as above and 
identify all fibers over $\partial D$
with $(F,B)$ so that $\sigma$ is constant. 
The set of all such
choices of identification (agreeing with the given one on $F_{[1:0]}$ and 
$TF|B$, up to fiberwise symplectic isotopy fixing $TF|B$)
is then given by $\pi_1(\D_{\omega_F})$. 
Passing to $\pi_1(\D)$ classifies Lefschetz pencils $f$ as in
Proposition \ref{prop:fixedA}, and $\sigma$ has a constant extension over
the remaining fibers of any $f$ coming from $\Im i_* \subset
\pi_1(\D)$ (for any choice of element in the corresponding
coset of $\ker i_*$).

Next we apply Theorem \ref{thm:taming}.
By contractibility of the space of $\sigma$--tame complex structures
on each $T_x F_{f(x)}$, we obtain a $\sigma$--tame family of
almost-complex structures on the fibers.
After declaring a suitable horizontal distribution to be complex,
we obtain a fiberwise $\sigma$--tame complex
structure $J$ on $X$, which we may assume agrees with $J_V$ near $B$
and with the structures on the chosen standard charts near the critical points.
Now $(f,J)$ is a linear 1--system as required.
Set $\nu= (\ker d\pi)|B$.
For each $F_y$, define $\eta_y$ on a neighborhood $W_y$ by pulling
back $\sigma|F_y$ by a map $r\co W_y\to W_y$ collapsing $W_y$ onto $F_y$ 
away from $B\cup K$ (cf proof of Theorem \ref{thm:taming}). 
Then each $\eta_y$ agrees with $\omega_V$ on a fixed
neighborhood of $B$, and with $\omega_j$ on a
neighborhood of the critical point if $F_y$ is singular.
We can assume each $\eta_y |D_j\cap W_y =0$.
Let $\zeta$ be any form on $X$  representing $c_f$, agreeing with
$\eta_{[1:0]}$ near $F_{[1:0]}$, and vanishing on each thimble.
(Note $c_f |F_{[1:0]} = [\omega_F] = [\eta_{[1:0]}]|F_{[1:0]}$ as
required, and the thimbles add no 2--homology to $W_{[1:0]}$ since $n\geq 3$.)
The condition $[\eta_y-\zeta] = 0 \in H_{\dR}^2 (W_y,B)$ is trivially
true for $y=$[1:0].
The case of any regular value $y$ then follows since $F_y$ comes with
an isotopy $\rel B$ in $X$ to $F_{[1:0]}$, sending $\eta_y|F_y$
to $\omega_F$.
For a critical value, we can assume $W_y$ is obtained from a tubular
neighborhood of a regular fiber by adding an $n$--handle.
Since $n\ge3$, the handle adds no 2--homology, so the condition holds
for all $y$.
Now Theorem~\ref{thm:taming} and Addendum~\ref{Jfixed} provide a unique
isotopy class of symplectic forms $\omega$ on $X$ taming $J$,
with $[\omega]=c_f$ and $\omega |F_{[1:0]}$ pairwise isotopic to $\omega_F$.
This $\omega$ can be assumed to satisfy the required conditions 
for Theorem~\ref{thm:omega}:
Compatibility of $\omega$ with $J^* =J$ is given on $\nu$.
It follows on $K$ since $\omega$ is made
from $\omega_t = t\eta + f^*\omega_{\text{std}}$; the second term
vanishes on $K$, and the first 
agrees with each $t\omega_j$ if we set $\{\hat y_1,\ldots,\hat y_m\} = f(K) 
\cup \{\text{[1:0]}\}$ when applying Theorem~\ref{thm:smoothmap}. 
Similarly, the thimbles $D_j$ are Lagrangian, since
$f^*\omega_{\text{std}}$ vanishes on them, as does $\eta$ if the forms
$\alpha_y$ arising from Theorem~\ref{thm:smoothmap} are chosen to vanish there.
(This can be arranged since $H^1_{\dR}(S^{n-1})=0$ for $n\geq 3$,
and $\alpha_{[1:0]} =0$.)
The forms $\widetilde \omega_s = \psi_s^* \eta_{y_0}$ 
in the deformation constructed for (1) also restrict to scalar
multiples of $\eta=0$ on each $S_j$ as required.

For fixed $X$ and $f$, we wish to compare two arbitrary
symplectic structures $\omega_u$, $u=0,1$, satisfying the conclusions
of the theorem.
We adapt the previous procedure to 1--parameter families, beginning
with the construction of $\sigma$.
For $u=0,1$,
choose $\pi_u\co V_u \to B$ with fibers tangent to the $\omega_u$--normal
bundle $\nu_u \to B$, and construct structures $\omega_{V, u}$ and
$J_{V, u}$ as before, using the Hermitian form $\omega_u|\nu_u$.
By Weinstein (cf proof of Theorem~\ref{thm:taming}),
we can assume $\omega_u = \omega_{V,u}$ near $B$ 
after an arbitrarily $C^1$--small isotopy.
(First isotope $\pi_u$ to get equality on $F_{[1:0]}$, preserving the 
fibers of $f$, then isotope $\omega$ fixing $F_{[1:0]}$.)
Let $\sigma_u$ be the family obtained by restricting $\omega_u$ to the
fibers.
Symplectic parallel transport gives a fiber-preserving map $\varphi_u\co
A \times F_{[1:0]} \to X$ for which $\varphi_u^*\sigma_u$ is constant,
and the Lagrangian thimbles $D_j$ are horizontal.
Now smoothly extend $\pi_u$, $\omega_{V,u}$, $J_{V,u}$ and $\varphi_u$
for $0 \le u \le 1$.
Also extend $\omega_u$ near $K$ by linear interpolation, so that it is
$J^*$--compatible on $K$ and has Lagrangian disks $\Delta_j$, $0\le u \le 1$.
Condition (1) gives a 
pairwise isotopy from $\omega_0$ through $\omega_F$ to $\omega_1$ on
$(F_{[1:0]}, B)$. 
Use this to extend $\sigma_u$ for $0 \leq u \leq 1$ to
$F_{[1:0]}$  with $\sigma_{1/2}=\omega_F$, 
and then extend as before to the fibers over $A$ and $D$, agreeing with
$\omega_{V,u}$ near $B$ and $\omega_u$ near $K$, and with $\varphi_u^*
\sigma_u$ constant for each $u$.

To complete the construction of $\sigma_u$, $u\in I=[0,1]$, over the fibers 
of $\id_I\times f$ on $I\times X$, we attempt to fill the hole 
over $(0,1)\times (\CP^1  -D)$, encountering obstructions.
As before, we can smoothly identify all fibers over $I\times \partial
D$ with $F_{[1:0]}$ so that the family $\sigma_u$ is constant for each $u$.
To fix this identification $\tau_u$ 
for each $u$, pull the preassigned path $\gamma$
back from $\D_{\omega_F}$ to $\D_{\sigma_u|F_{[1:0]}}$ by the given isotopy.
This moves the spheres $S_j$, but the same isotopy shows how to restore
them to their original position through families of Lagrangian spheres
(by the last part of (1)), yielding a canonically induced path $\gamma_u$
from the required representative of $\prod t_j$ to $\delta_{\sigma_u}$ in
$\D_{\sigma_u|F_{[1:0]}}$, which determines $\tau_u$.
Now let $\widetilde\tau_0$ be the identification of all fibers over 
$\{0\}\times (\CP^1 - \Int D)$ with $F_{[1:0]}$ by $\omega_0$--symplectic 
parallel transport along straight lines to [1:0] in $\CP^1-\Int D$. 
Comparing $\widetilde\tau_0|\partial D$ with $\tau_0$, 
we obtain the element of $\pi_1(\D)$ classifying
$f$, and see that this must lie in $\Im i_*$ (being explicitly represented 
by a loop in $\D_{\sigma_0|F_{[1:0]}} \simeq \D_{\omega_F}$). 
The corresponding construction of $\widetilde\tau_1$ from $\omega_1$ 
also gives $f$, so  $\widetilde\tau_0|\partial D$ and $\widetilde\tau_1|\partial D$
differ by an element $\beta$ of $\ker i_*$ 
(after we extend each $\widetilde\tau_u$ over $I\times \partial D$ to 
$u=1/2$ so that we can work with the fixed symplectic form $\omega_F$).
Similarly, a direct comparison of $\widetilde\tau_0$ and $\widetilde\tau_1$ 
yields an element 
$\alpha\in\pi_2 (\D,\D_{\omega_F}) \cong \pi_2(\D/\D_{\omega_F})$ with 
$\partial_*\alpha =\beta$. 
If $\alpha=0$, we can extend $\sigma_u$ over $X$ for $0\le u\le1$.
If only $\beta=0$, then we can perturb $\widetilde\tau_1$ so that 
$\alpha\in \pi_2(\D)$, and $\alpha$ provides a self-diffeomorphism of $X$
preserving $f$ and fixing $f^{-1}(D)$, after which $\sigma_u$ extends. 
These vanishing conditions are also necessary for the deformation 
and symplectomorphism, respectively, specified by the theorem, since
any allowable deformation $\omega_u$, $0\le u\le 1$, 
determines a family $\widetilde\tau_u$ as above interpolating between 
$\widetilde\tau_0$ and $\widetilde\tau_1$, showing that $\alpha=0$. 
(Note that the family $\omega_u$ comes with a smooth
family of isotopies from $\omega_u|F_{[1:0]}$ to $\omega_F$ as in (1),
allowing us to continuously pull back $\gamma$ to each 
$\D_{\sigma_u|F_{[1:0]}}$ as before, to define the required interpolation 
$\tau_u$ between $\tau_0$ and $\tau_1$ in the absence of the condition 
$\sigma_{1/2}|F_{[1:0]} = \omega_F$.) 

To complete the proof for $n\ge3$, it suffices to  
construct the required deformation between $\omega_0$ and $\omega_1$
from the completed family $\sigma_u$ on $I\times X$.
First find a continuous family $J_u$ of fiberwise $\sigma_u$--tame
almost-complex structures on $X$ as before, using a
horizontal distribution that is $\omega_u$--orthogonal to the fibers
when $u=0,1$.
Then $J_u$ is $\omega_u$--tame, $u=0,1$.
For $0< u< 1$, the previous argument now produces  suitable symplectic
structures $\omega_u$ on $X$, with $\sigma_u$ replacing $\omega_F$ in (1).
The family $\omega_u$, $0\le u\le1$, need not be continuous (particularly
at 0,1), but we can smooth it as for Theorem~\ref{thm:taming},  
with a partition of unity on $I$,
obtaining the required deformation $\omega'_u$:
First pull each $\omega_u$ back to a neighborhood $I_u$ of $u\in I$,
by a map preserving the bundles $\nu_u$, forms $\omega_B$ and disks $D_j$.
For $I_u$ sufficiently small, the required conditions (1--4) are preserved,
where the isotopy in (1) from $\omega_u$ to $\sigma_u$ on $F_{[1:0]}$
at each $v\in I_u$ is through $J_{v}$--taming forms
(by Addendum~\ref{Jfixed}) and is defined to be constant for $u=0,1$.
Splicing by a partition of unity on $I$ preserves (2--4), producing
the required family $\omega'_u$ and a 2--parameter family of
symplectic structures on $F_{[1:0]}$, interpolating between $\omega'_u$
and a convex combination $\sigma'_u$ of nearby forms
$\sigma_{v}$ for each $u$.
If the intervals $I_u$ were sufficiently small, we can extend by
$(1-s)\sigma'_u + s\sigma_u$ to obtain a  2--parameter family of
symplectic structures from $\omega'_u$ to $\sigma_u$, all
agreeing with $\omega_F$ on $B$ and each $S_j$, and constant for $u=0,1$.
Moser's technique, parametrized by $u$, now produces a 2--parameter family
of diffeomorphisms of $(F,B)$, which can be reinterpreted as a smooth
family of isotopies as in (1)
from $\omega'_u|F_{[1:0]}$ to $\omega_F$ $=\sigma_{1/2}|F_{[1:0]}$,
interpolating between the given ones for $\omega_0,\omega_1$.

For the remaining case $n\le 2$, inclusion $\D_{\omega_F}\subset\D$ 
is a homotopy equivalence, so we must show each $f$ has a unique deformation
class of structures $\omega$ as specified.
The case $n=1$ ($f\co X\to \CP^1$ a simple branched covering) is
trivial, so we assume $n=2$.
Then $f$ is a hyperpencil,
so \cite[Theorem 2.11(b)]{G1} gives the required form $\omega$, provided
we arrange $J^*$--compatibility as before.
Uniqueness of the deformation class follows the method of proof of 
\cite[Theorem 1.4]{G2} with $m=0$ and $K$ replaced by $K\cup B$:
Given two forms $\omega_u$, $u=0,1$, as specified in Theorem~3.3, 
we can find an $(\omega_{\text{std}},f)$--tame, $\omega_u$--tame almost-complex 
structure $J_u$ for each $u$ (cf also \cite[Lemma 2.10]{G1}).
We interpolate to a family $J_u$, $0\le u\le 1$ 
(eg by contractibility in \cite[Theorem 2.11(a)]{G1}),
and construct the required deformation $\omega_u$, $0\le u\le 1$, using 
a partition of unity on $I$ as before.
Conditions (1--4) for each $\omega_u$ are easily verified, with only 
(1) requiring comment:  We are given 
isotopies $\rel B$ from $\omega_u |F_{[1:0]}$ to $\omega_F$, $u= 0,1$. 
Moser provides an isotopy $\rel B$ from $\omega_0$ to $\omega_1$ on 
$F_{[1:0]}$. 
Combining these three  isotopies gives a path representing an element of
$\pi_1(\D',\D'_{\omega_F})$, where the prime indicates rotations of
$TF|B$ are allowed.
This group vanishes, allowing us to extend our path into the required 
2--parameter family of diffeomorphisms.
\end{proof}

\section*{Acknowledgement}
The author was partially supported by NSF grant DMS--0102922.

\Addresses\recd

\end{document}